\documentstyle[twoside,12pt]{article}
\newtheorem{prop}{Proposition}[section]
\newtheorem{dfn}[prop]{Definition}
\newtheorem{theo}[prop]{Theorem}
\newtheorem{conj}[prop]{Conjecture}
\newtheorem{rem}[prop]{Remark}
\newtheorem{coro}[prop]{Corollary}

\newtheorem{exam}[prop]{Example}

\title{Dual Cones and Mirror Symmetry for
Generalized Calabi-Yau Manifolds}

\author{Victor V. Batyrev$\,^*$ \\
 Universit\"at-GHS-Essen, Fachbereich  6,  Mathematik \\
 Universit\"atsstr. 3,  45141  Essen  \\
Federal Republic of Germany \\
e-mail: matf$\emptyset \emptyset$@vm.hrz.uni-essen.de \\
and \\
Lev A. Borisov \\
Department of Mathematics \\
University of Michigan \\
Ann Arbor, Michigan 48109-1003 \\
e-mail: Lev.Borisov@math.lsa.umich.edu}

\date{}

\begin{document}

\maketitle

\footnote{$\,^*$ Supported by DFG, Forschungsschwerpunkt  Komplexe
Mannigfaltigkeiten.}

\begin{abstract}
We introduce a special class of convex rational polyhedral cones which
allows to construct generalized Calabi-Yau varieties of
dimension $(d + 2(r-1))$, where $r$ is a positive integer and
$d$ is the dimension of critical string vacua with cenral charge
$c = 3d$. It is conjectured that the natural combinatorial
duality satisfied by  these  cones
corresponds to the  mirror involution. Using the theory of toric
varieties, we show that our conjecture includes
as special cases all already known examples of mirror pairs
proposed  by physicists and agrees with previous conjectures
of the authors concerning explicit constructions
of  mirror manifolds. In particular we obtain a mathematical
framework which explains the construction of
mirrors of rigid Calabi-Yau manifolds.
\end{abstract}

\newpage

\section{Introduction}
\noindent

This paper is devoted to the problem of finding an appropriate
mathematical framework which explains explicit constructions
of mirror pairs of Calabi-Yau manifolds. In this context, the existence
of rigid Calabi-Yau manifolds furnish a fundamental obstruction to
a mathematical formulation of mirror symmetry motivated
by a natural involution in $N=2$  superconformal field theories.
In recent works \cite{schimm1,schimm2}, Schimmrigk  came
to the conclusion that the class of Calabi-Yau manifolds
{\em is not appropriate setting for mirror symmetry}.
For every positive integer $r$,  Schimmrigk has proposed
a new class of K\"ahler manifolds of dimension $(d + 2(r-1))$
which generalizes the class of $d$-dimensional Calabi-Yau varieties.
Although the first Chern class of
{\em generalized Calabi-Yau manifolds} does not vanish for $r >1$, it is
possible to derive from these manifold the massless spectrum of
critical string vacua. Typical examples of the generalized Calabi-Yau
manifolds are obtained from quasi-smooth  hypersurfaces
of degree $w$ in weighted projective spaces
\[ {\bf P}(w_1, \ldots, w_{d + 2r}) \]
where the degree $w$ and the weights $w_1, \ldots, w_{d + 2r}$  are
related by the condition:
\[ w_1 +  \cdots + w_{d + 2r} = rw. \]
Schimmrigk has shown that some classes of usual $d$-dimensional
Calabi-Yau manifolds $V$ which can be described
as complete intersections of $r$ hypersurfaces in products of
$r$ copies of weighted projective spaces naturally give rise to
generalized Calabi-Yau manifolds $M$ of dimension
$(d- 2(r-1))$ embedded as hypersurface
in some higher-dimensional weighted projective space. Moreover,
the correspondence $V \mapsto M$ induces the canonical
inculison of the Hodge $(p,q)$-spaces
\[ H^{p,q}(V) \subset H^{p +r -1, q + r -1}(M). \]
It is important to remark that not every generalized Calabi-Yau manifold
of dimension $(d + 2(r-1))$
can be obtained by this method from some usual $d$-dimensional
Calabi-Yau manifold.

\medskip

In \cite{bat.dual,bat.straten,borisov} there were proposed
a combinatorial approach to the construction
of mirror pairs using theory of toric varieties \cite{danilov,oda}.
The main idea of
this approach is the interpretation of the mirror duality in terms of the
classical duality for convex sets. It is a natural ask whether the same
approach can be applied to rigid and generalized Calabi-Yau manifolds?
The main purpose of this paper is to show that the answer is positive.
\bigskip

In Section 2, we define a special class of convex rational polyhedral
cones which we call {\em reflexive Gorenstein cones}.
Every reflexive Gorenstein cone $\sigma$ canonically defines the
projective toric Fano variety
${\bf P}_{\sigma}$ together with the ample invertible sheaf
${\cal O}_{{\bf P}_{\sigma}}(1)$ on it such that
some $r$-tensor power of ${\cal O}_{{\bf P}_{\sigma}}(1)$ is
isomorphic to the anticanonical sheaf of ${\bf P}_{\sigma}$ (in particular
${\bf P}_{\sigma}$  has only Gorenstein singularities).
The zeros of global sections of ${\cal O}_{{\bf P}_{\sigma}}(1)$ are
generalized Calabi-Yau manifolds.  For special simplicial
reflexive Gorenstein cones, we obtain generalized Calabi-Yau
manifolds considered by Schimmrigk.
The class of reflexive Gorenstein cones of fixed dimension admits a
natural involution $\sigma \rightarrow \check{\sigma}$ which we conjecture
to correspond the mirror involution in $N =2$ superconformal theories.
The rest of the paper is devoted to arguments which confirm this
conjecture.

In Section 3, we give a general overview for
the method of reduction of complete
intersections in toric varieties to hypersurfaces in higher dimensional
toric varieties.  This method allows to construct  reflexive
Gorenstein cones and
generalized Calabi-Yau manifolds from usual Calabi-Yau complete
intersections in Gorenstein toric Fano varieties. The relation
between Hodge structures of Calabi-Yau complete intersections
and the corresponding generalized Calabi-Yau manifolds enable
to use results from \cite{bat.var}
to determine  variations of Hodge structure for Calabi-Yau complete
intersections in toric varieties.

In Section 4, we check that the duality between nef-partitions
defining Calabi-Yau complete intersections in Gorenstein toric
Fano varieties \cite{borisov} agree with the duality between
reflexive Gorenstein cones.

In Section 5, we show that the duality for reflexive Gorenstein cones
agree also with the explicit construction of mirrors for rigid Calabi-Yau
manifolds.
\medskip

{\bf Acknowledgements.} We are grateful to I.V. Dolgachev, H. Esnault,
S. Katz, Yu.I. Manin, D. van Straten, and R. Schimmrigk for helpful
discussions.

\bigskip

\section{Reflexive Gorenstein cones}
\noindent

Let $\overline{M}$ and $\overline{N} =
{\rm Hom}(\overline{M}, {\bf Z})$ be dual free abelian groups
of rank $\overline{d}$, $\overline{M}_{\bf R}$ and
$\overline{N}_{\bf R}$ the real scalar extensions
of $\overline{M}$ and $\overline{N}$,
$ \langle *, * \rangle\; : \; \overline{M}_{\bf R}
\times \overline{N}_{\bf R} \rightarrow
{\bf R} $
the natural pairing.
We consider $\overline{M}$ (resp. $\overline{N}$)  as
a maximal lattice in $\overline{M}_{\bf R}$
(resp. in $\overline{N}_{\bf R}$).

\begin{dfn} {\rm
Let $e_1, \ldots, e_k$ be elements of $M$.
By a finite rational polyhedral cone
\[\sigma = {\bf R}_{\geq 0} \langle e_1, \ldots, e_k \rangle
\subset \overline{M}_{\bf R} \]
generated by $\{ e_1, \ldots, e_k \}$ we mean the  set of all $x =
\lambda_1 e_1 + \cdots + \lambda_k e_k \in \overline{M}_{\bf R}$ such that
$\lambda_i \geq 0$, $\lambda_i \in {\bf R}$ $(i =1, \ldots, k)$.
}
\end{dfn}

\begin{dfn}
{\rm Let $\sigma$ be a
finite rational polyhedral cone in $\overline{M}_{\bf R}$.
Then  the {\em dual cone} is  the set
\[ \check{\sigma} = \{ y \in \overline{N}_{\bf R} \mid
\langle x, y \rangle \geq 0 \mbox{ for all } x \in \sigma \}. \] }
\end{dfn}

\begin{dfn} {\rm
If a  $\overline{d}$-dimensional cone
$\sigma \subset \overline{M}_{\bf R}$ satisfies
the condition $\sigma \cap -\sigma
=0$, then we can  uniquely choose
the  {\em minimal set of integral generators} of $\sigma$ which
is defined as the set of all primitive $\overline{M}$-lattice vectors
on $1$-dimensional faces of $\sigma$. By a {\em generator} of
$\sigma$ we will always mean a primitive $\overline{M}$-lattice
vector on a $1$-dimensional face one of $\sigma$.  }
\end{dfn}

\begin{dfn}
{\rm Let $\sigma$ be a finite rational polyhedral cone in
$\overline{M}_{\bf R}$ satisfying the condition
$\sigma \cap -\sigma =0$.
Then the  cone $\sigma$ is called {\em Gorenstein} if
there exists an element $n_{\sigma}  \in N$ such that
$\langle e, n_{\sigma} \rangle =1$ for each generator $e$ of $\sigma$.}
\end{dfn}

\begin{rem}
{\rm If $\sigma \subset \overline{M}_{\bf R}$ is a Gorenstein cone
whose dimension equals the dimension of $\overline{M}_{\bf R}$,
then the element $n_{\sigma}$ is uniquely  defined. }
\end{rem}

\begin{dfn}
{\rm A Gorenstein cone $\sigma \subset \overline{M}_{\bf R}$
is called {\em reflexive} if the
dual cone $\check{\sigma} \subset \overline{N}_{\bf R}$
is also Gorenstein. In this case,
the positive integer
\[ r_{\sigma} = \langle m_{\check{\sigma}}, n_{\sigma} \rangle \]
will be called the {\em index} of $\sigma$ (or of $\check{\sigma})$.
Since the notion of reflexive Gorenstein cone strongly depends on
the choice of the maximal sublattice $\overline{M}$ in the
$\overline{d}$-dimensional real vector space, we will say
that  reflexive Gorenstein cones are defined by pairs
$(\sigma, \overline{M})$. }
\end{dfn}
\medskip

The notions of Gorenstein cones and reflexive Gorenstein cones
can be interpreted via toric geometry.

\begin{prop}
Let ${\bf A}_{\sigma}= {\rm Spec}\,{\bf C} \lbrack \check{\sigma} \cap
\overline{N} \rbrack$ be the affine $\overline{d}$-dimensional
toric variety associated with a rational polyhedral
cone $\sigma$.  Then $\sigma$ is a Gorenstein cone if and only if
${\bf A}_{\sigma}$ has only Gorenstein singularities.
\label{gorenstein}
\end{prop}

{\em Proof. } The statement follows from the characterization of
M. Reid for Gorenstein toric singularities \cite{reid}. \hfill $\Box$
\medskip

By a lattice polyhedron in a finite dimensional real vector space
${\cal U}$
over ${\bf R}$  we always mean a convex polyhedron whose vertices
belong to some fixed maximal sublattice in ${\cal U}$.

\begin{dfn}
{\rm Let $\sigma$ be a Gorenstein cone. Then
\[ \Delta_{\sigma} = \{ x \in \sigma \mid \langle x, n_{\sigma}
\rangle =1 \} \]
is a $(\overline{d} -1)$-dimensional convex lattice polyhedron
which we call the {\em support} of $\sigma$. }
\end{dfn}

\begin{rem} {\rm
Every lattice polyhedron $\Delta$ in a $d$-dimensional
space $M_{\bf R}$ can be considered as a support of
$(d+1)$-dimensional Gorenstein cone $\sigma_{\Delta} \subset
\overline{M}_{\bf R} = {\bf R} \oplus M_{\bf R}$
defined as
\[ \sigma_{\Delta} = \{
( \lambda, \lambda x) \in \overline{M}_{\bf R}
\mid \lambda \in {\bf R}_{\geq 0}, \; x \in \Delta \}. \]
}
\end{rem}

\begin{dfn}
{\rm \cite{bat.dual} A lattice polyhedron is called {\em reflexive} if
$\sigma_{\Delta}$ is a reflexive Gorenstein cone of
index $1$. If $\Delta$ is a reflexive polyhedron, then
the support of the dual reflexive Gorenstein cone $\check{\sigma}_{\Delta}$
is another reflexive polyhedron $\Delta^*$ which is called {\em dual} to
$\Delta$. }
\label{def.ref}
\end{dfn}

\begin{prop}
A Gorenstein cone $\sigma$ is reflexive cone of index $r$
if and only $r\Delta_{\sigma}$ is a reflexive polyhedron.
\label{reflex.cone}
\end{prop}

{\em Proof. } Let $\sigma \subset \overline{M}_{\bf R}$ be a
Gorenstein cone and $r$ a positive
integer. Define new lattices $M' \subset \overline{M}_{\bf R}$ and
$N' \subset \overline{N}_{\bf R}$ as follows
\[ M' = \{ x \in \overline{M} \mid
\langle x, n_{\sigma} \rangle = 0\,(\mbox{\rm mod
$r$ }) \},  \]
\[ N' = \overline{N} + {\bf Z} n_{\sigma}'\; \mbox{\rm where }\; n_{\sigma}'
= \frac{1}{r}n_{\sigma}. \]
Then the lattice $N'$ is dual to the lattice
$M'$ and the pair $(\sigma, M')$ defines again a Gorenstein cone
whose support is $r\Delta_{\sigma}$.

Assume now that the pair $(\check{\sigma}, \overline{N})$ defines a reflexive
Gorenstein cone of index $r$.  Since $m_{\check{\sigma}} \in M'$ and
$\langle m_{\check{\sigma}}, e \rangle = 1$ for every $\overline{N}$-integral
generator $e$ of ${\check{\sigma}}$, we obtain that $e$ is also
a $N'$-integral generator of ${\check{\sigma}}$.
Thus  the pair $(\check{\sigma}, N')$ defines a reflexive
Gorenstein cone of index $1= \langle m_{\check{\sigma}}, n_{\sigma}'\rangle$.
 By \ref{def.ref}, $r\Delta_{\sigma}$ is reflexive.
Analogous arguments show the part "if".
\hfill $\Box$

 The next statement follows immediately from the
equivalent characterizations of reflexive polyhedra in \cite{bat.dual}:

\begin{prop}
A lattice polyhedron $\Delta$ is reflexive if and only if
$\Delta$ is the support of global sections of the ample anticanonical
sheaf on a Gorenstein toric Fano variety.
\label{reflex}
\end{prop}

Using  \ref{reflex.cone} and \ref{reflex}, we obtain:

\begin{coro} Let $\sigma$ be a Gorenstein cone.
We define the degree of an element
$m \in \sigma \cap \overline{M}$ as
${\rm deg}\, m = \langle m, n_{\sigma} \rangle$.
Let ${\bf P}_{\sigma}= {\rm Proj}\,{\bf C} \lbrack {\sigma} \cap
\overline{M} \rbrack$ be the corresponding
projective $(\overline{d}-1)$-dimensional
toric variety. Then $\sigma$ is reflexive Gorenstein cone
of index $r$ if and only if
${\bf P}_{\sigma}$ is a Gorenstein toric Fano variety such that
the anticanonical sheaf on ${\bf P}_{\sigma}$ is isomorphic to
${\cal O}_{{\bf P}_{\sigma}}(r) =
{\cal O}_{{\bf P}_{\sigma}}(1)^{\otimes r}$.
\label{characteriz}
\end{coro}

By the adjunction formula, one has:

\begin{coro}
Let $\sigma$ be a reflexive Gorenstein cone. Then the zeros of
a general global section of ${\cal O}_{{\bf P}_{\sigma}}(1)$ define
a $(\overline{d} -2)$-dimensional algebraic variety $Z$ with
only Gorenstein toroidal singularities such that
the anticanonical sheaf on $Z$ is ${\cal O}_Z(r-1)$.
\label{zeros}
\end{coro}

\begin{dfn}
{\rm We call the $(\overline{d} -2)$-dimensional algebraic varieties
$Z$ obtained as zeros of global sections of ${\cal O}_{{\bf P}_{\sigma}}(1)$
(\ref{zeros}) {\em generalized Calabi-Yau manifolds associated
with the reflexive Gorenstein cone $\sigma$.}}
\end{dfn}

\begin{exam}
{\rm Let ${\bf P}(w_1, \ldots, w_{\overline{d}})$ be a
$(\overline{d}-1)$-dimensional weighted projective spaces whose
weights satisfy the condition $w_1 + \cdots + w_{\overline{d}} = r w_0$
for some positive integers $r$ and $w_0$ such that $w_i$ divides
$w_0$ for all $i =1, \ldots, \overline{d}$.
Let $\sigma = {\bf R}^d_{\geq 0} \subset {\bf R}^{\overline{d}}$
be the positive octant.
Define the maximal lattice $\overline{M} \subset {\bf R}^{\overline{d}}$
as
\[ \overline{M} = \{ (x_1, \ldots, x_{\overline{d}}) \in
{\bf Z}^{\overline{d}} \mid w_1 x_1 + \cdots + w_{\overline{d}}
x_{\overline{d}} = 0\, (\mbox{\rm mod } w_0) \}. \]
Then the pair $(\sigma, \overline{M})$ defines an example
of a reflexive Gorenstein cone of index $r$.
In this case the associated with $\sigma$
generalized Calabi-Yau manifolds
are hypersurfaces of degree $w_0$ in
${\bf P}(w_1, \ldots, w_{\overline{d}})$. }
\label{weighted}
\end{exam}

Now we formulate the main conjecture:

\begin{conj}
Every pair of $\overline{d}$-dimensional
 dual refelxive Gorenstein cones $\sigma$ and $\check{\sigma}$
of index $r$ give rises  to  a $N=2$ superconformal theory with central
charge $c = 3( \overline{d} - 2(r-1))$. Moreover, the superpotentials
of the corresponding Landau-Ginzburg  theories define two
families of generalized Calabi-Yau manifolds
associated with $\sigma$ and $\check{\sigma}$ which are exchanged by
the mirror involution.
\label{conjecture}
\end{conj}

\begin{rem}
{\rm By \ref{def.ref}, the duality between reflexive Gorenstein
cones is equivalent to the duality  between the supporting reflexive polyhedra.
This shows that Conjecture \ref{conjecture}
includes the one of \cite{bat.dual}.  }
\end{rem}
\bigskip

\section{From complete intersections to hypersurfaces}

\noindent

In this section we discuss the general procedure which ascribes to a complete
 intersection in toric variety a hypersurface in another toric variety.
 This method is essentially due to Danilov and Khovanski\v{i} \cite{dan.khov}
 (see also \cite{esnault,terasoma}).
Then we show that the hypersurfaces that arise in this way give rise
to reflexive Gorenstein cones if and only if the original complete intersection
 is Calabi-Yau variety (Prop.  \ref{reflexivecone}).
We also discuss when the the reflexive Gorenstein cones come from
some Calabi-Yau complete intersections.
\medskip

Let $M$ be a free abelian group of rank $d$, $M_{\bf R}$ the real scalar
extension of $M$.
Let $\Delta_1, \ldots, \Delta_r \subset M_{\bf R}$ be
lattice polyhedra  supporting global sections of
semi-ample invertible sheaves
${\cal O}_X(D_1), \ldots, {\cal O}_X(D_r)$
on a $d$-dimensional toric variety $X$.
Without loss of generality we will always assume that ${\rm dim}\,
\Delta_1 + \cdots + \Delta_r = d$.
Let ${\bf Z}^r$ be the standard  $r$-dimensional lattice,
${\bf R}^r$ its real scalar
extension.  We put $\overline{M} = {\bf Z}^r \oplus
M$, $\overline{d} = d + r$,  and define
the cone $\sigma \subset \overline{M}_{\bf R}$
as
\[ \sigma = \{ ( \lambda_1, \ldots, \lambda_r, \lambda_1 x_1 + \cdots
+ \lambda_r x_r ) \in \overline{M}_{\bf R}
\mid \lambda_i \in {\bf R}_{\geq 0}, \; x_i \in \Delta_i, \;
i =1, \ldots r \}. \]

 Then $\sigma$ is a $\overline{d}$-dimensional
Gorenstein cone, where $n_{\sigma}$ is an element of the dual lattice
$\overline{N}$ defined uniquely by the conditions
\[ \langle x, n_{\sigma} \rangle = 0 \; \mbox{ for } x \in M_{\bf R} \subset
\overline{M}_{\bf R}; \]
\[ \langle e_i, n_{\sigma} \rangle = 1\; \mbox{ for } i = 1, \ldots, r, \]
where $\{e_1, \ldots, e_r \} $ is the standard basis of ${\bf Z}^r \subset
\overline{M}$.

\begin{rem}
{\rm The supporting polyhedron $\Delta_{\sigma}$ coincides
with the $(\overline{d} -1)$-dimensional polyhedron
$\Delta_1 * \cdots * \Delta_r$ considered by
Danilov and Khovanski\v{i} (\cite{dan.khov} \S 6). }
\end{rem}

\begin{prop}
Denote by $Y$ the $(\overline{d} -1)$-dimensional toric variety which is
the toric ${\bf P}^{r-1}$-bundle over $X:$
\[ Y = {\bf P}({\cal O}_X(D_1) \oplus \cdots \oplus {\cal O}_X(D_r)). \]
Let  ${\cal O}_Y(-1)$ be  the Grothendieck tautological sheaf on $Y$.
Then ${\cal O}_Y(1)$ is semi-ample, and
${\bf P}_{\sigma}$ is the birational image of the toric  morphism
\[ \alpha \; :\; Y \rightarrow {\bf P}_{\sigma} \]
defined by global sections of ${\cal O}_Y(1)$.
In particular the polyhedron $\Delta_{\sigma}$ supports the
global sections of ${\cal O}_Y(1)$.
\label{reduction}
\end{prop}

{\em Proof.} Let $\pi\, :\, Y \rightarrow X$ be the
canonical projection. Since $\pi$ agrees with the
torus actions on $X$ and $Y$, we obtain the natural torus action
on ${\cal O}_Y(1)$. Since
\[ \pi_* {\cal O}_Y(1) =  {\cal O}_X(D_1) \oplus \cdots
\oplus {\cal O}_X(D_r)  \]
is the direct sum of sheaves generated by global sections,
${\cal O}_Y(1)$ is also generated by global sections; i.e.,
${\cal O}_Y(1)$ is semi-ample.
In order to determine the polyhedron supporting the global sections
of ${\cal O}_Y(1)$, it suffices to compute the $(d+r-1)$-dimensional
torus action on
\[ H^0({\cal O}_Y(1)) \cong H^0({\cal O}_X(D_1) \oplus \cdots
\oplus {\cal O}_X(D_r)). \]
The latter is trivial, since we know  the $(r-1)$-dimensional
torus action on homogeneous coordinates in ${\bf P}^{r-1}$
and the $d$-dimensional
torus action on $H^0({\cal O}_X(D_i))$ $(i =1, \ldots, r)$
defined by the lattice points in the polyhedron $\Delta_i$.
\hfill $\Box$

\begin{coro}
Every global section $s$ of  ${\cal O}_Y(1)$ defines uniquely
global sections $s_i$ of ${\cal O}_X(D_i)$ such that
\[ \pi_*( s) = (s_1, \ldots , s_r) \in
H^0({\cal O}_X(D_1)) \oplus \cdots \oplus H^0({\cal O}_X(D_r)), \]
and vise versa, every $r$-tuple of sections $s_i \in H^0({\cal O}_X(D_i))$
$( i =1, \ldots, r)$ defines a global section $s \in
H^0({\cal O}_Y(1))$.
\label{sections}
\end{coro}

\begin{coro}
Let $s \in H^0({\cal O}_Y(1))$ and
$s_i \in H^0({\cal O}_X(D_i))$
$( i =1, \ldots, r)$ be global sections as in \ref{sections}.
Denote by $V_s$ the hypersurface in $Y$ defined by $s = 0$, and by
$V_{s_i}$ $( i =1, \ldots, r)$ the hypersurfaces  in $Y$ defined by
$s_i = 0$. Then $Y \setminus V_s$ is locally trivial in Zariski topology
${\bf C}^{r-1}$-bundle over
\[ X \setminus \bigcap_{i=1}^r V_{f_i}. \]
\label{bundle}
\end{coro}

\begin{rem}
{\rm The statement in \ref{bundle} implies the isomorphism
\[ H^{i}_c( X \setminus \bigcap_{i=1}^r V_{f_i}) \cong
H^{i + 2r -2}_c (Y \setminus V_s) \]
which sends the $(p,q)$-component in
cohomology with compact supports of
$X \setminus \bigcap_{i=1}^r V_{f_i}$
to $(p + r -1,q+ r-1)$-component in the cohomology with
compact supports of $Y \setminus V_s$.
The relation between Hodge structures of the complements to the
higher dimensional hypersurfaces and to the complete
intersections was used   for
estimations of the Hodge type of algebraic subvarieties \cite{esnault}
and  in the proof of weak global Torelli theorem  \cite{terasoma}.
One can  consider \ref{bundle} also  as a version of the Lagrange method
proposed by Danilov and Khovanski\v{i} (\cite{dan.khov}, \S 6) for
affine hypersurfaces}.
\end{rem}

\begin{prop}
The  cone $\sigma \subset \overline{M}_{\bf R}$ is a
reflexive Gorenstein cone of index $r$ if and only if
${\cal O}_X(D_1 + \cdots + D_r)$
is isomorphic to the anticanonical sheaf
on $X$.
\label{reflexivecone}
\end{prop}

{\em Proof.} Standard calculations show that the
canonical sheaf ${\cal K}_Y$ on $Y$ is isomorphic to the tensor product
\[ {\cal O}_Y(-r) \otimes \pi^*  {\cal K}_X \otimes
\pi^* \Lambda^r({\cal O}_X(D_1) \oplus \cdots
\oplus {\cal O}_X(D_r)). \]
Since
\[ \Lambda^r({\cal O}_X(D_1) \oplus \cdots
\oplus {\cal O}_X(D_r)) \cong {\cal O}_X(D_1 + \cdots + D_r), \]
we have
\[ {\cal K}_Y^{-1} \cong {\cal O}_Y(r) \otimes \pi^*  {\cal K}_X^{-1}
\otimes \pi^* {\cal O}_X(-D_1 - \cdots - D_r). \]

Assume that ${\cal O}_X(D_1 + \cdots + D_r)$ is isomorphic to
the anticanonical sheaf ${\cal K}^{-1}_X$
on $X$. Then ${\cal O}_Y(r)$ is isomorphic to the anticanonical sheaf on
$Y$. Therefore, by \ref{characteriz} and \ref{reduction},
$\sigma$ is reflexive.

The "only if" part is left to reader. \hfill $\Box$
\medskip

The following  example of the reduction of Calabi-Yau
complete intersections to a higher-dimensional generalized
Calabi-Yau manifolds is due to Schimmrigk \cite{schimm2}.

\begin{exam}
{\rm Let $k, l$ be two positive integers. Define $V$ as a
Calabi-Yau complete intersection of two hypersurfaces in ${\bf P}^k
\times {\bf P}^l$ having bidegrees $(k+1, 1)$ and $(0,l+1)$. Then
the corresponding generalized Calabi-Yau manifolds associated
with the reflexive Gorenstein cone of index $2$ are hypersurfaces
of degree  $(k+1)l$ in the $(k+l-1)$-dimensional weighted projective space
\[ {\bf P} (\underbrace{(l-1),\ldots,(l-1)}_{k+1},
\underbrace{(k+1), \ldots,(k+1)}_{l+1}).\]
}
\end{exam}

\begin{dfn}
{\rm A $\overline{d}$-dimensional  reflexive Gorenstein cone
$\sigma \subset \overline{M}_{\bf R}$  is
called {\em  split} if there exist two lattice polyhedra
$\Delta_1, \Delta_2 \subset M_{\bf R}$ of
dimension $d = \overline{d} - 2$
such that
\[ \sigma \cong \{ ( \lambda_1,  \lambda_2, \lambda_1 x_1 +
\lambda_2 x_2 ) \in \overline{M}_{\bf R}
\mid \lambda_i \in {\bf R}_{\geq 0}, \; x_i \in \Delta_i, \;
i =1, 2 \}. \]}
\label{split1}
\end{dfn}

\begin{dfn}
{\rm A $\overline{d}$-dimensional  reflexive Gorenstein cone
$\sigma \subset \overline{M}_{\bf R}$ of index $r$ is
called {\em completely split} if there exist $r$ lattice polyhedra
$\Delta_1, \ldots, \Delta_r \subset M_{\bf R}$ of
dimension $d = \overline{d} - r$
such that
\[ \sigma \cong \{ ( \lambda_1, \ldots, \lambda_r, \lambda_1 x_1 + \cdots
+ \lambda_r x_r ) \in \overline{M}_{\bf R}
\mid \lambda_i \in {\bf R}_{\geq 0}, \; x_i \in \Delta_i, \;
i =1, \ldots r \}. \]}
\label{split2}
\end{dfn}

\begin{rem}
{\rm Cones that come from Calabi-Yau
complete intersections are completely split. On
 the other hand, every splitting of the reflexive Gorenstein cone gives rise to
 the family of complete intersections.
However it's not clear whether there exist reflexive Gorenstein cones that
could
 be split in several essentially different ways, so that they come from
 different Calabi-Yau varieties.}
 \label{diffCYsamecones}
 \end{rem}

There exist simple examples of reflexive Gorenstein cones of
index $r >1$ which are not split (and hence are not completely split):

\begin{exam}
{\rm Let $\sigma = {\bf R}_{\geq 0}^{\overline{d}} \subset
{\bf R}^{\overline{d}}$ be the positive octant. Assume that
$\overline{d} = kr$ where $k, r$ are positive integers and
$k,r >1$.  Define the lattice $\overline{M}$ as
\[ \overline{M} = \{ (x_1, \ldots, x_{\overline{d}}) \in
{\bf Z}^{\overline{d}} \mid x_1 + \cdots + x_{\overline{d}} = 0\,
\mbox{\rm mod $k$} \}. \]
Then the pair $(\sigma, \overline{M})$ defines a reflexive
Gorenstein cone of index $r$ which is not split. Indeed,
if there were two $(\overline{d} - 2)$-dimensional polyhedra
$\Delta_1, \Delta_2$ having  the property described in
\ref{split1}, then for any two vertices $v_1 \in \Delta$,
$v_2 \in \Delta_2$ we could find two generators
$e_1 = (1,0,v_1)$, $e_2 = (0,1,v_2)$ of  the cone $\sigma$ such
that the segment $\lbrack e_1, e_2 \rbrack$ would have no
interior $\overline{M}$-lattice points. On the other hand,
it is clear that for any two generators of $\sigma$
the segment $\lbrack e_1, e_2 \rbrack$ always
contains $k-1$ interior $\overline{M}$-lattice points. }
\end{exam}

\begin{rem}
{\rm In Section 5 we consider an example of a $3d$-dimensional
cone reflexive Gorenstein cone
$\sigma$ which is completely split, but the dual cone
$\check{\sigma}$ is not split.}
\end{rem}

\section{Complete intersections and nef-partitions}

\noindent

There is an important class of reflexive Gorenstein cones $\sigma$
such that both $\sigma$ and $\check{\sigma}$ are completely split.
These cones correspond to so called {\em nef-partitions} introduced
in \cite{borisov}.

We use notations from  the
previous section and assume that $X$ is a Gorenstein
 Fano toric variety. Let $T \subset X$ be the dense open torus
orbit, $E_1, \ldots, E_k$ irreducible components of
$X \setminus T$. Then ${\cal O}(E_1 + \cdots + E_k)$ is
is  naturally isomorphic to the anticanonical sheaf ${\cal K}_X^{-1}$
on $X$. Denote by $I$ the set $\{ 1, \ldots, k \}$. Put
$E=\sum_{i \in I} E_i.$

\begin{dfn}
{\rm The decomposition of the index set $I$ into a disjoint union of $r$ sets
 $I_j,\,j=1, \ldots, r$ is called a {\em nef-partition} if all
  $D_j=\sum_{i\in I_j}E_i$ are
 semi-ample Cartier divisors on $X$. By abuse of notations, we
 will call by nef-partition also the set of convex lattice polyhedra
 $\Pi = \{ \Delta_1, \ldots, \Delta_r \}$ such that $\Delta_i$ is
 the support of global sections of ${\cal O}(D_i)$ $i =1, \ldots, r$. }
\label{defnef}
\end{dfn}

\begin{dfn}
{\rm Let $\Pi = \{ \Delta_1, \ldots, \Delta_r \}$ be a nef-partition,
$\Sigma$ is the fan in $N_{\bf R}$ defining the Gorenstein toric
variety $X$. For every $E_i$ ($ i =1, \ldots, k$), we denote by $e_i$
the primitive $N$-lattice generator of the $1$-dimensional cone of
$\Sigma$ corresponding to the divisor $E_i$. Define the lattice
polyhedron $\nabla_j$ $(j =1, \ldots, r)$ as the convex hull
\[ \nabla_i = {\rm Conv}(\{0\} \cup \bigcup_{j \in J_i} \{ e_j \}). \]
}
\label{defdual}
\end{dfn}

\begin{rem}
{\rm Since $D_i = \sum_{j \in J_i} E_i$ defines a convex piecewise linear
function $\psi_i$ such that $\psi_i(e_j) = 1$ if $j \in J_i$ and
$\psi_i(e_j) = 0$ otherwise, we can define $\Delta_1,
\ldots, \Delta_r$ as
\[ \Delta_i = \{ x \in M_{\bf R} \mid \langle x, y \rangle \geq -
\psi_i(y) \}, \; i=1, \ldots, r. \]
In the sequel, we will always assume this definition which
 immediately implies that all polyedra
$\Delta_1, \ldots, \Delta_r$ contain $0 \in M_{\bf R}$. }
\end{rem}

The main result of \cite{borisov} is the following:

\begin{theo}
The set $\Pi^* = \{ \nabla_1, \ldots, \nabla_r \}$ is also a nef-partition.
In other words, there exist  another $d$-dimensional Gorenstein toric Fano
variety $X^*$ which compactifies the dual torus $T^*$ and the index set
$J^* = \{ 1, \ldots , l \}$ for the irreducible components
$E_1^*,  \ldots,  E_l^*$ of  $X^* \setminus T^*$
 such that the $\nabla_i$ is the support of global sections
 of the semi-ample sheaf ${\cal O}(D_i^*) =
 {\cal O} (\sum_{j \in J^*_i} E_i^*)$ $( i =1, \ldots, r)$ where
 $J_1^* \cup \cdots \cup J_r^*$ is a splitting of $I^*$ into a
 disjoint union.
 \end{theo}

 The combinatorial involution on $\Pi \mapsto \Pi^*$  is conjectured
 to give rise to the mirror symmetry for the families of Calabi-Yau complete
 intersections in $X$ and $X^*$. Some results which  confirm this conjecture
 are obtained in \cite{bat.bor}.

 On the other hand, once we have
\[ {\cal K}^{-1}_X = {\cal O}_X(D_1 + \cdots + D_r) \]
 and
\[ {\cal K}^{-1}_{X^*} = {\cal O}_{X^*}(D_1^* + \cdots + D_r^*) \]
we can follow the pro\-ce\-du\-re of the pre\-vious sec\-tion
and get two ref\-le\-xive $\overline{d}$-di\-men\-sio\-nal
Gorenstein cone $\sigma \subset \overline{M}_{\bf R}$ and
$\sigma^* \subset \overline{N}_{\bf R}$.

The purpose of this section is to relate the involution for
nef-partitions  \cite{borisov} to the involution
for  reflexive Gorenstein cones.

First  we recall one property
 of dual nef-partitions:

\begin{prop} {\rm \cite{borisov}}
Let $x \in \Delta_i$, $y \in \nabla_j$. Then
\[ \langle x , y \rangle \geq -1\;\;   \mbox{ if $ i = j$}, \]
\[ \langle x , y \rangle \geq 0\;\;   \mbox{ if $ i \neq j$}. \]
\label{property}
\end{prop}

The main statement is contained in the following theorem.

\begin{theo}
 Assume that the canonical pairing $\langle \cdot , \cdot \rangle \, :
M \times N \rightarrow {\bf Z}$ is extended to the pairing
between $\overline{M} = {\bf Z}^r \oplus M$  and $\overline{N} = {\bf Z}^r
\oplus N$ as
\[ \langle (a_1, \ldots,a_r, m), (b_1, \ldots, b_r, n) \rangle =
\sum_{i =1}^r a_i b_i + \langle m, n \rangle. \]
Then the cone $\sigma$ is dual to $\sigma^*$; i.e., the
dual nef-partitions in the sense of {\rm \cite{borisov}} give rise to
the dual reflexive Gorenstein cones.
\label{nefpart-refcone}
\end{theo}

{\em Proof.} By \ref{property}, if $(a_1, \ldots,a_r, m) \in \sigma$ and
$(b_1, \ldots, b_r, n) \in \sigma^*$, then
\[ \langle (a_1, \ldots,a_r, m), (b_1, \ldots, b_r, n) \rangle \geq 0. \]
Therefore $\sigma^* \subset \check{\sigma}$.

Let $(b_1, \ldots, b_r, n) \in \check{\sigma}$. Then, for any
$x \in \Delta_i$,  one has $\langle x, n \rangle \geq - b_i$
$(i = 1, \ldots, r)$. Using the same arguments as in the proof
of Prop. 3.2 in \cite{borisov}, we obtain that
$ n \in b_1 \nabla_1 + \cdots + b_r \nabla_r$. Therefore $\check{\sigma}
\subset \sigma^*$.
 \hfill $\Box$

\begin{rem}
{\rm Several different nef-partitions can give rise to the same reflexive
Gorenstein cone, so the duality for nef-partitions carries more
information than that for the cones. From the geometrical point of view,
the nef-partition corresponds to the family of Calabi-Yau
complete intersections
 together with some special degeneration of the
family into the union of the strata of dim $d-r$. It is not yet clear
if the dual family depends upon the degeneration or not, which is related to
the
 question in Remark \ref{diffCYsamecones}. Of course, there are no such
 problems in the case of hypersurfaces, or, equivalently, when the index of the
 reflexive Gorenstein cones is $1$.}
\end{rem}

\section{Mirrors of rigid Calabi-Yau manifolds}
\noindent

All already known constructions  of  the mirror correspondence
for rigid Calabi-Yau manifolds were originated from
 the explicit identification of
the minimal superconformal models of Gepner \cite{gepner} with
the special Landau-Ginzburg superpotentials \cite{greene.vafa.warner}. This
indentification allows to apply the orbifolding  modulo some finite
symmetry group \cite{greene.plesser,roan}.
Our purpose  is to  show  that the orbifold-construction of mirrors
of rigid Calabi-Yau varieties agrees with the duality for reflexive
Gorenstein cones.
\bigskip

We consider in details the example
of the rigid $d$-dimensional Calabi-Yau manifold associated with
the superconformal theory $1^{3d}$ which is the tensor
product of $3d$ copies of the level-1 theories.
\medskip

Let $E_0 ={\bf C} / {\bf Z} \langle 1, \tau \rangle$
be the unique elliptic curve having an authomorphism of
order $3$; i.e., $J(E_0) = 0$, $\tau = e^{\pi i/3}$
and $E_0$ is isomorphic to the Fermat cubic in ${\bf P}^2$.
Notice that
the action of the group ${\bf Z}/3{\bf Z}$ on $E_0$ has exactly $3$
fixed points which we denote by $p_0, p_1, p_2$. Let
\[ G = \{ (g_1,\ldots, g_d) \in ({\bf Z}/3{\bf Z})^d \mid
g_1 + \cdots g_d = 0\, (\mbox{mod } d) \}. \]
Then $G$ is the maximal subgroup in $({\bf Z}/3{\bf Z})^d$
whose action on the product $X = (E_0)^d$
leaves invariant the holomorphic $d$-form
$z_1 \wedge \cdots \wedge z_d$. We denote by
$Z$ the geometric quotient $X/G$ considered as a $d$-dimensional
orbifold.
Then the mirror involution in the $1^{3d}$ superconformal
theory shows that $Z$ is mirror symmetric to the $(3d-2)$-dimensional
Fermat cubic $Y \subset {\bf P}^{3d-1}$. One sees a
geometric confirmation of this  duality from
the following  statement:

\begin{prop}
Let $\hat{Z}$ be a maximal projective crepant partial resolution of
quotient singularities of $Z$ {\rm \cite{bat.dual}}. Then
\[ h^{1,1} (\hat{Z}) = h^{3d-3,1}(Y). \]
\end{prop}

{\em Proof.} Using the standard technique based on the consideration
of the Jacobian ring associated with the homogeneous equation of
$Y$, we obtain that a basis of $H^{3d-3,1}(Y)$ can be identified with
all square-free monomials of degree $3$ in $3d$ variables. Therefore
\[ h^{3d-3,1}(Y) = { 3d \choose 3 } = \frac{ d (3d-1)(3d-2) }{2}. \]
On the other hand, there exists exactly $d$ linearly independent
$G$-invariant $(1,1)$-forms on $X$: $dz_1\wedge d \overline{z}_1,
\ldots, dz_d \wedge d\overline{z}_d$. Thus, we have $h^{1,1}(Z) = d$. It
remains to compute the number of exceptional divisors on $\hat{Z}$.

Let $\gamma_1 \,: \, \hat{Z} \rightarrow Z$ be the maximal projective
crepant partial resolution, $\gamma_2\,: \,  X \rightarrow Z$ the
quotient by $G$.
One easily sees that the $\gamma_1$-image of an exceptional
divisor $D \subset \hat{Z}$
in $Z$ is either $(d-2)$-dimensional, or $(d-3)$-dimensional.
In the first case, $\gamma_1(D)$ is the $\gamma_2$-image of
a $G$-invariant codimension-2 subvariety defined by
the conditions $z_i, z_j \in \{ p_0, p_1, p_2\}$
$(i \neq j)$, and the $\gamma_1$-fiber
over general point of $\gamma_1(D)$ is the exceptional locus
of the crepant resolution of the $2$-dimensional
Hirzebruch-Jung singularity of type
$A_2$, i.e., it consists of two irreducible components. In the second
 case, $\gamma_1(D)$ is the $\gamma_2$-image of
a $G$-invariant codimension-3 subvariety defined by
conditions $z_i, z_j, z_k \in \{ p_0, p_1, p_2\}$,
and the $\gamma_1$-fiber
over general point of $\gamma_1(D)$
is an irreducible surface.
Therefore the number of the exceptional  divisors equals
\[ 2 \cdot 3^2 \cdot { d \choose 2} + 3^3 \cdot { d \choose 3 }. \]
This immediately implies
\[ h^{1,1}(\hat{Z}) =  d + 2 \cdot 3^2 \cdot
{ d \choose 2} + 3^3 \cdot { d \choose 3 } = \frac{ d (3d-1)(3d-2) }{2}. \]
\hfill $\Box$
\medskip

The combinatorial interpretation of the mirror duality between
$Y$ and $Z$ is based on the representation of $X$ and $Z$ as
complete intersections in toric varieties. Let ${\bf P}_{\Delta}$
be the $2$-dimesnional toric variety associated with
the reflexive polygon $\Delta = {\rm Conv}\{ (1,0),(0,1),(-1,-1) \}$.
We can also define ${\bf P}_{\Delta}$ in ${\bf P}^3$
by the equation $u_0^3 = u_1 u_2 u_3$.

\begin{prop}
Let $C \subset {\bf P}_{\Delta}$ be a curve defined by an
equation $\lambda_1 u_1 + \lambda_2 u_2 + \lambda_3 u_3 = 0$.
Then $C$ is isomorphic to $E_0$.
\end{prop}

{\em Proof. } It is sufficient to notice that the mapping
\[  q \;: \; {\bf P}_{\Delta}  \rightarrow {\bf P}_{\Delta} \]
\[ (u_0, u_1,u_2,u_3) \mapsto (e^{2\pi i/3}u_0, u_1,u_2,u_3) \]
induces an authomorphism of order $3$ of $C$ with three fixed points. \hfill
 $\Box$

\begin{coro}
The $d$-dimensional variety $X$ is a complete intersection
of $d$ nef-divisors in the $2d$-dimensional
toric variety $({\bf P}_{\Delta})^d$.
\end{coro}

\begin{coro}
The mapping $q$ induces the action of $({\bf Z}/3{\bf Z})^d$ on
$({\bf P}_{\Delta})^d$ such that $Z$ becomes a
complete intersection in the geometric quotient
$({\bf P}_{\Delta})^d/G$.
\end{coro}

\begin{prop}
Let $\sigma \subset {\bf R}^{3d}$ be positive octant. Define
\[ \overline{M} = {\bf Z}^{3d} + {\bf Z}(\frac{1}{3}, \ldots,
\frac{1}{3}). \]
Then the pair $(\sigma, M)$ defines a
$3d$-dimensional Gorenstein reflexive cone associated with
$Z$ as a complete intersection in $({\bf P}_{\Delta})^d/G$.
\end{prop}

{\em Proof. } We notice that the $3$-dimensional
reflexive cone $\sigma_{\Delta}$ can be described
as the positive octant in ${\bf R}^3$ with respect to
the lattice
\[ M = {\bf Z}^3 + {\bf Z}(\frac{1}{3},\frac{1}{3},\frac{1}{3}). \]
Thus the pair $((\sigma_{\Delta})^d, M^d)$ is the
$3d$-dimensional Gorenstein reflexive cone associated with
$X$ as a complete intersection in $({\bf P}_{\Delta})^d$.
It remains to compute the sublattice $\overline{M} \subset M^d$
corresponding to modding out of $({\bf P}_{\Delta})^d$ by $G$.
It is clear that ${\bf Z}^{3d} \subset \overline{M}$ and
$M^d/\overline{M}$ must be isomorphic to
$G$. On the other hand, by construction,
$\overline{M}$ must be invariant under
$d$-element permutations in $M^d$. These conditions define
$\overline{M}$ uniquely as
${\bf Z}^{3d} + {\bf Z}({1}/{3}, \ldots, {1}/{3})$. \hfill
$\Box$.

Using \ref{weighted}, we have:

\begin{coro}
Let ${\sigma} \subset {\bf R}^{3d}$
be positive octant. Define
\[ \overline{N} = \{ (x_1, \ldots, x_{3d}) \in
{\bf Z}^{3d} \mid x_1 + \cdots + x_{3d} = 0( \mbox{\rm mod 3}). \]
Then the pair $(\sigma, N)$ defines the dual to $(\sigma, M)$
reflexive Gorenstein cone with respect to the standard scalar
product on ${\bf R}^{3d}$. In particular,
the reflexive Gorenstein pair $( \sigma, \overline{N})$
corresponds to cubic hypersurfaces in ${\bf P}^{3d-1}$.
\end{coro}

\begin{rem}
{\rm If we choose an intermediate lattice $M'\, : \;
\overline{M} \subset M' \subset M^d$ and the corresponding dual intermediate
lattice $N'\, : \; N^d \subset N' \subset \overline{N}$, then we obtain
another pair of dual reflexive Gorenstein cones. In particular,
if $d =3$, then one  can obtain a rigid Calabi-Yau 3-fold $Z'$ by modding out
$E_0 \times E_0 \times E_0$ by the diagonal action of ${\bf Z}/3{\bf Z}
\subset G \cong ({\bf Z}/3{\bf Z})^2$. The
mirror symmetric generalized Calabi-Yau manifolds are then obtained
as quotients of $7$-dimensional cubics by ${\bf Z}/3{\bf Z}$. This particular
case was  considered in \cite{cand.derr.parkes}.}
\end{rem}

\end{document}